\documentclass{amsart}
\newtheorem{Thm}{Theorem}[section]
\newtheorem{Prop}[Thm]{Proposition}
\newtheorem{Lem}[Thm]{Lemma}
\newtheorem{Cor}[Thm]{Corollary}
\newtheorem{Ex}[Thm]{Example}
\newtheorem{Prob}[Thm]{Problem}
\theoremstyle{remark}
\newtheorem{Rem}[Thm]{Remark}
\numberwithin{equation}{section}

\DeclareMathOperator{\supp}{supp}
\begin{document}

\title[Zero divisors and $L^p(G)$, II]
{Zero divisors and $L^p(G)$, II}

\author[P. A. Linnell]{Peter A. Linnell}
\address{Math \\ VPI \\ Blacksburg \\ VA 24061--0123
\\ USA}
\email{linnell@math.vt.edu}
\urladdr{http://www.math.vt.edu/people/linnell/}

\author[M. J. Puls]{Michael J. Puls}
\address{New Jersey City University \\ Jersey City \\
NJ 07305--1597 \\ USA}
\email{mpuls@njcu.edu}
\urladdr{http://ellserver3.njcu.edu/math/puls/Puls.htm}

\begin{abstract}
Let $G$ be a discrete group, let $p \ge 1$,
and let $L^p(G)$ denote the
Banach space $\{\sum_{g\in G} a_g g
\mid \sum_{g\in G} |a_g|^p < \infty\}$.  The following problem will
be studied: given $0 \ne \alpha \in \mathbb {C}G$ and
$0 \ne \beta \in L^p(G)$, is $\alpha * \beta \ne 0$?  We will
concentrate on the case $G$ is a free abelian or
free group.
\end{abstract}

\keywords{zero divisor, free group, Fourier transform, radial
function, free abelian group}

\subjclass{Primary: 43A15; Secondary: 43A25, 42B99}

\date{March 27, 2000}
\maketitle

\section{Introduction} \label{Sintroduction}

Let $G$ be a discrete group and let $f$ be a complex-valued function
on $G$.  We may represent $f$ as a formal sum $\sum_{g \in G} a_g g$
where $a_g \in \mathbb {C}$ and $f(g) = a_g$.  Thus $L^{\infty}(G)$
will consist of all formal sums $\sum_{g\in G} a_g g$ such that
$\sup_{g \in G} |a_g|  < \infty$, $C_0(G)$ will consist of those formal
sums for which the set $\{ g \mid |a_g| > \epsilon \}$ is finite for
all $\epsilon > 0$, and for $p \ge 1$,
$L^p(G)$ will consist of those formal sums for which $\sum_{g \in G}
|a_g|^p < \infty$.  Then we have the following inclusions:
\[
\mathbb {C}G \subseteq L^p(G) \subseteq C_0(G) \subseteq L^{\infty}
(G).
\]
For $\alpha = \sum_{g\in G} a_g g \in L^1(G)$ and $\beta =
\sum_{g\in G} b_g g \in L^p(G)$, we define a multiplication
$L^1(G) \times L^p(G) \to L^p(G)$ by
\begin{equation} \label{Emult}
\alpha * \beta = \sum_{g,h} a_g b_h gh =
\sum_{g\in G}\left( \sum_{h\in G} a_{gh^{-1}}b_h \right) g.
\end{equation}
In this paper we consider the following:
\begin{Prob}
Let $G$ be a torsion free group and let $1 \le p \le \infty$.
If $0 \ne \alpha \in \mathbb {C}G$ and $0 \ne \beta \in
L^p(G)$, is $\alpha * \beta \ne 0$?
\end{Prob}
Some results on this problem are given in \cite{linnell,
puls}.  In this sequel we shall obtain new results for the
cases $G = \mathbb {Z}^d$, the free abelian group of rank $d$,
and $G= F_k$, the free group of rank $k$.

Part of this work was carried out while
the first author was at the
Sonderforschungsbereich in M\"unster.  He would like to thank
Wolfgang L\"uck for organizing his visit to M\"unster, and the
Sonderforschungsbereich for financial support.

\section{Statement of Main Results}

Let $0 \ne \alpha \in L^1(G)$ and let $1 \le p \in \mathbb {R}$.
We shall say that
$\alpha$ is a \emph{$p$-zero divisor}
if there exists $\beta \in L^p(G)
\setminus 0$ such that $\alpha * \beta = 0$.  If $\alpha * \beta
\ne 0$ for all $\beta \in C_0(G) \setminus 0$, then we say that
$\alpha$ is a \emph{uniform nonzero divisor}.

Let $2 \le d \in \mathbb {Z}$.  It was shown in \cite{puls}
that there are $p$-zero divisors in $\mathbb {C}\mathbb {Z}^d$
for $p > \frac{2d}{d-1}$.  In this paper we shall show that
this is the best possible by proving
\begin{Thm} \label{Tabelian1}
Let $2 \le d \in \mathbb {Z}$, $1 \le p \in \mathbb {R}$,
let $0 \ne \alpha \in \mathbb {C}
\mathbb {Z}^d$, and let $0 \ne \beta \in L^p(\mathbb {Z}^d)$.
If $p \le \frac{2d}{d-1}$, then $\alpha * \beta \ne 0$.
\end{Thm}
Let $\mathbb {T}^d$ denote the $d$-torus which, except in
Section \ref{Skey}, we will
view as the cube $[-\pi, \pi]^d$ in $\mathbb {R}^d$ with opposite
faces identified, and let $\mathfrak {p} \colon [-\pi, \pi]^d \to
\mathbb {T}^d$ denote the natural surjection.
For $n \in \mathbb {Z}^d$ and $t \in \mathbb {T}^d$,
let $n \cdot t$ indicate the dot product, which is well defined
modulo $2 \pi$.  If $\alpha = \sum_{n \in \mathbb {Z}^d}
a_n n \in L^1 (\mathbb {Z}^d)$, then for $t \in \mathbb {T}^d$
its Fourier transform
$\hat{\alpha} \colon \mathbb {T}^d \to \mathbb {C}$
is defined by
\[
\hat{\alpha} (t) = \sum_{n \in \mathbb {Z}^d} a_ne^{-i(n \cdot t)}
\]
and we shall let $Z (\alpha) = \{ t \in \mathbb {T}^d \mid
\hat{\alpha} (t) = 0 \}$.  We say that $M$ is a hyperplane in
$\mathbb {T}^d$ if there exists a hyperplane $N$ in $\mathbb{R}^d$
such that $M = \mathfrak {p} ([-\pi, \pi]^d \cap N)$.
We will prove the following theorem,
which is an improvement over \cite[theorem 1]{puls}.
\begin{Thm} \label{Tabelian2}
Let $\alpha \in \mathbb {C} \mathbb {Z}^d$.
Then $\alpha$ is a uniform nonzero divisor if and only if
$Z (\alpha)$ is
contained in a finite union of hyperplanes in $\mathbb {T}^d$.
\end{Thm}

Let $V = \mathfrak{p} \bigl( (-\pi, \pi)^d \bigr)$, let $\alpha \in
L^1(\mathbb {Z}^d)$, let $E = Z (\alpha) \cap V$, and let $U$ be
an open subset of $(-\pi, \pi)^{d-1}$.  Let $\phi \colon U \to (-\pi,
\pi)$ be a smooth map, and suppose $\{ \mathfrak {p}(x, \phi (x))
\mid x \in U \} \subseteq E$.  If the Hessian matrix
\[
\left( \frac{\partial^2 \phi}{\partial x_i \partial x_j} \right)
\]
of $\phi$ has constant rank $d-1-\nu$ on $U$ where
$0 \le \nu \le d-1$,
then we say that $\phi$ has constant relative nullity $\nu$.  We
shall say that $Z (\alpha)$ has \emph{constant relative nullity
$\nu$} if every localization $\phi$ of $E$ has constant relative
nullity $\nu$ \cite[p.~64]{guo}.  We shall prove
\begin{Thm}  \label{Tabelian3}
Let $\alpha \in \mathbb {C} \mathbb {Z}^d$, let $1 \le p \in \mathbb
{R}$, and let $2 \le d \in
\mathbb {Z}$.  Suppose that $Z (\alpha)$ is a
smooth $(d-1)$-dimensional submanifold of $\mathbb {T}^d$
with constant relative nullity $\nu$ such that $0 \le \nu \le d-2$.
Then $\alpha$ is a $p$-zero divisor if and only if $p > \frac{2
(d-\nu)}{d-1-\nu}$.
\end{Thm}

For $k \in \mathbb {Z}_{\ge 0}$, let
$F_k$ denote the free group on $k$ generators.  It was proven in
\cite{linnell} that if $\alpha \in \mathbb {C}F_k \setminus 0$ and
$\beta \in L^2(F_k) \setminus 0$, then $\alpha * \beta \ne 0$.  We will
give an explicit example to show that if $k \ge 2$, then this result
cannot be extended to $L^p(F_k)$ for any $p > 2$.  This is a bit
surprising in view of Theorem \ref{Tabelian1}.  We will conclude this
paper with some results about $p$-zero divisors for the free group
case.

\section{A characterization of $p$-zero divisors}

Let $G$ be a group, not necessarily discrete, and let $L^p(G)$ be the
space of $p$-integrable functions on $G$ with respect to Haar
measure.  Let $y \in G$ and let $f \in L^p(G)$.  The right
translate of $f$ by $y$ will be denoted by $f_y$, where $f_y(x) =
f(xy^{-1})$.  Define $T^p[f]$ to be the closure in $L^p(G)$ of all
linear combinations of right translates of $f$.  A common problem
is to determine when $T^p[f] = L^p(G)$; see
\cite{edgarrosenblatt,edwards,rosenblatt} for background.

Now suppose that $G$ is also discrete.  Given
$1 \le p \in \mathbb {R}$,
we shall always let $q$ denote the conjugate index of $p$.  Thus
if $p > 1$, then $\frac {1}{p} + \frac {1}{q} = 1$, and if $p=1$ then
$q=\infty$.  Sometimes we shall require $p= \infty$, and then $q = 1$.
Let $\alpha = \sum_{g \in G} a_g g \in L^p(G)$, $\beta =
\sum_{g \in G} b_gg \in L^q(G)$, and define a map
$\langle \cdot ,\cdot
\rangle \colon L^p(G) \times L^q(G) \to \mathbb {C}$ by
\[
\langle \alpha, \beta \rangle = \sum_{g \in G} a_g \overline{b_g}.
\]
Fix $h \in L^q(G)$.  Then $\langle \cdot ,h \rangle$ is a continuous
linear functional on $L^p(G)$ and if $p \ne \infty$, then every
continuous linear functional on $L^p(G)$ is of this form.
We shall use the notation
$\tilde {\beta }$ for  $\sum_{g \in G} b_g g^{-1}$,
$\bar {\beta }$ for  $\sum_{g \in G} \overline {b_g} g$, and
$\beta^*$ for $\sum_{g\in G} \overline{b_g} g^{-1}$.
Also the same formula in equation \eqref{Emult} gives a
multiplication $L^p(G) \times L^q(G) \to L^{\infty}(G)$.
Then we have the following elementary lemma,
which roughly says that $\alpha * \beta = 0$ if and only if all the
translates of $\alpha$ are perpendicular to $\beta$.

\begin{Lem} \label{Lperp}
Let $1 \le p \in \mathbb {R}$ or $p = \infty$,
let $\alpha \in L^p(G)$, and let
$\beta \in L^q(G)$.  Then $\alpha * \beta = 0$ if and only if
$\langle (\tilde {\alpha})_y, \bar{\beta} \rangle = 0$ for all $y \in G$.
\end{Lem}
\begin{proof}
Write $\alpha = \sum_{g\in G} a_g g$ and $\beta = \sum_{g\in G} b_g
g$.  Then
\[
\alpha * \beta = \sum_{y \in G} (\sum_{g \in G} a_{y g^{-1}} b_g)y
\]
and $\langle (\tilde{\alpha})_y, \bar{\beta}
\rangle = \sum_{g \in
G} a_{yg^{-1}} b_g$.  The result follows.
\end{proof}

The following proposition, which is a
generalization of \cite[lemma 1]{puls}, characterizes $p$-zero
divisors in terms of their right translates (the statement of
\cite[lemma 1]{puls} should have the additional condition that $p \ne
1$).

\begin{Prop} \label{Pchar}
Let $\alpha \in L^1(G)$ and let $1 < p \in \mathbb {R}$ or
$p=\infty$.  Then $\alpha$ is a $p$-zero divisor if and
only if $T^q[\tilde{\alpha}] \ne L^q(G)$.
\end{Prop}
\begin{proof}
The Hahn-Banach theorem tells us that $T^q[\tilde {\alpha}]
\ne L^q(G)$ if and
only if there exists a nonzero continuous linear functional on
$L^q(G)$ which vanishes on $T^q[\tilde {\alpha}]$.  The result
now follows from Lemma \ref{Lperp}.
\end{proof}

\begin{Rem}
If $p=1$ in the above Proposition \ref{Pchar}, we would need to
replace $L^q(G)$ with $C_0(G)$, and $T^q[\tilde{\alpha}]$ with the
closure in $C_0(G)$ of all linear combinations of right translates of
$\tilde{\alpha}$.
\end{Rem}

\section{A key proposition} \label{Skey}

In this section we prove a proposition that will enable us to prove
Theorems \ref{Tabelian1}, \ref{Tabelian2} and \ref{Tabelian3}.

Let $1 \le p \in \mathbb {R}$, let
$y \in \mathbb {R}^d$ and let $f \in L^p(\mathbb {R}^d)$.  We
shall use additive notation for the group operation in $\mathbb
{R}^d$; thus the right translate of $f$ by $y$ is now given by $f_y
= f(x-y)$.  We say that $f$ has linearly independent translates if
and only if for all $a_1, \dots ,a_m \in \mathbb {C}$, not all zero,
and for all distinct $y_1, \dots , y_m \in \mathbb {R}^d$,
\[
\sum_{i=1}^m a_if_{y_i} \ne 0.
\]
For the rest of this section we shall view $\mathbb {T}^d$ as the
unit cube $[0,1]^d$ with opposite faces identified.  Let $L^p(\mathbb
{T}^d \times \mathbb {Z}^d)$ denote the space of functions on
$\mathbb {T}^d \times \mathbb {Z}^d$ which satisfy
\[
\int_{t \in \mathbb {T}^d} \sum_{m \in \mathbb {Z}^d} |f(t,m)|^p \, dt
< \infty.
\]
Then for $\alpha = \sum_{n \in \mathbb {Z}^d} a_nn \in \mathbb
{C}\mathbb {Z}^d$ and $f \in L^p(\mathbb {T}^d \times \mathbb
{Z}^d)$, we define $\alpha f \in L^p (\mathbb {T}^d \times \mathbb
{Z}^d)$ by
\[
(\alpha f)(t,m) = \sum_{n \in \mathbb {Z}^d} a_n f(t,m-n),
\]
and this yields an action of $\mathbb {C} \mathbb {Z}^d$ on
$L^p(\mathbb {T}^d \times \mathbb {Z}^d)$.

\begin{Lem}  \label{Llinind}
Let $\alpha \in \mathbb {C} \mathbb {Z}^d$.  Then there exists $\beta
\in L^p(\mathbb {Z}^d) \setminus 0$ such that $\alpha * \beta = 0$ if
and only if there exists $f \in L^p(\mathbb {T}^d \times \mathbb
{Z}^d) \setminus 0$ such that $\alpha f = 0$.
\end{Lem}
\begin{proof}
Let $\beta \in L^p(\mathbb{Z}^d) \setminus 0$ such that $\alpha * \beta
= 0$ and define a nonzero function $f \in L^p(\mathbb {T}^d \times
\mathbb {Z}^d)$ by $f(t,m) = \beta (m)$. For $n \in \mathbb {Z}^d$, set
$b_n = \beta(n)$.  Then
\begin{equation} \label{Elinind}
\begin{split}
(\alpha f)(t,m) &= \sum_{n \in \mathbb {Z}^d} a_n f(t, m-n)
= \sum_{n \in \mathbb {Z}^d } a_n \beta (m-n) \\
&= \sum_{n \in \mathbb {Z}^d} a_n b_{m-n}
= (\alpha * \beta ) (m) = 0.
\end{split}
\end{equation}

Conversely suppose there exists $f \in L^p( \mathbb {T}^d \times
\mathbb {Z}^d) \setminus 0$ such that $\alpha f = 0$.  This means
that $(\alpha f)(t,n) = 0$ for all $n$, for all $t$ except on a set
$T_1 \subset \mathbb {T}^d$ of measure zero.  Also
$\sum_{n \in \mathbb {Z}^d} |f(t,n)|^p  < \infty$
for all $t$ except on a set $T_2
\subset \mathbb {T}^d$ of measure zero.  Since $f \ne 0$, we may
choose $s \in \mathbb {T}^d
\setminus (T_1 \cup T_2)$ such that $f(s,n) \ne 0$ for some $n$.  Now
define $\beta (n) = f(s,n)$.  Then $\beta \in L^p(\mathbb {Z}^d)
\setminus 0$ and the calculation in equation
\eqref{Elinind} shows that $\alpha * \beta = 0$.
\end{proof}

For $\alpha = \sum_{n \in \mathbb {Z}^d} a_n n \in \mathbb {C}\mathbb
{Z}^d$ and $f \in L^p(\mathbb {R}^d)$, we define $\alpha f \in
L^p(\mathbb {R}^d)$ by
\[
(\alpha f)(x) = \sum_{n \in \mathbb {Z}^d} a_n f(x-n).
\]
If $\alpha \ne 0$ and $\alpha f = 0$, then there is a dependency
among the right translates of $f$, i.e.\ $f$ does not have linearly
independent translates.  We are now ready to prove

\begin{Prop} \label{Plinind}
Let $\alpha \in \mathbb {C} \mathbb {Z}^d$.  Then
$\alpha $ is a $p$-zero divisor if and only if there exists $f \in
L^p(\mathbb {R}^d) \setminus 0$ such that $\alpha f = 0$.
\end{Prop}
\begin{proof}
Define a Banach space isomorphism
$\zeta \colon L^p(\mathbb {R}^d) \to L^p(\mathbb {T}^d \times \mathbb
{Z}^d)$ by $(\zeta f)(t,n) = f(t-n)$ for $f \in L^p(\mathbb {R}^d)$.  We
want to show that this isomorphism commutes with the action of
$\mathbb {C} \mathbb {Z}^d$.  Clearly it will be sufficient to show
that $\zeta$ commutes with the action of $\mathbb {Z}^d$.  If $m \in
\mathbb {Z}^d$, then
\begin{alignat*}{2}
m \bigl((\zeta f)(t,n) \bigr) &= m(f(t-n)) &&= f(t-n-m) \\
&= (mf)(t-n) &&= (\zeta (mf))(t,n).
\end{alignat*}
Thus the action of $\mathbb {C}\mathbb {Z}^d$ commutes with $\zeta$.  We
deduce that for $\alpha \in \mathbb {C}\mathbb {Z}^d$, there exists
$f \in L^p(\mathbb {R}^d) \setminus 0$ such that $\alpha f = 0$ if
and only if there exists $ f' \in L^p(\mathbb {T}^d \times
\mathbb {Z}^d) \setminus 0$ such that $\alpha  f' = 0$.  The proposition
now follows from Lemma \ref{Llinind}.
\end{proof}

\begin{Rem} \label{Rlinind}
Replacing $L^p(\mathbb {R}^d)$ by $C_0(\mathbb {R}^d)$ in the above
arguments, we can also show that $\alpha$ is a uniform nonzero
divisor if and only if $\alpha f \ne 0$ for all
$f \in C_0(\mathbb {R}^d) \setminus 0$.
\end{Rem}

\section{Proofs of theorems \ref{Tabelian1}, \ref{Tabelian2},
and \ref{Tabelian3}}

The proof of Theorem \ref{Tabelian1} is obtained by combining
\cite[theorem 3]{rosenblatt} with Proposition \ref{Plinind}.  The
proof of Theorem \ref{Tabelian2} is obtained by combining
\cite[theorem 2.12]{edgarrosenblatt} with Remark \ref{Rlinind}.

Before we prove Theorem \ref{Tabelian3}, we will need to define the
notion of a $q$-thin set.  See \cite{edwards} for
more information on this and other concepts used in this paragraph.
Let $G$ be a locally compact abelian group and let $X$ be its
character group.  Let $\beta \in L^{\infty}(G)$ and let $\hat
{\beta}$ indicate the generalized Fourier transform of $\beta$.  The
key reason for using the generalized Fourier transform is that for
$\alpha \in L^1(G)$, we have $\widehat {\alpha * \beta} =
\hat{\alpha} \hat{\beta}$ which tells us that $\alpha * \beta = 0$ if
and only if $\supp \hat{\beta} \subseteq Z(\alpha)$.  Let $E
\subseteq X$.  We shall say that $E$ is $q$-thin if $\beta \in C_0(G)
\cap L^p(G)$ and $\supp \hat{\beta} \subseteq E$ implies $\beta =
0$.  Recall that $p$ is the conjugate index of $q$.  The result of
Edwards \cite[theorem 2.2]{edwards} says that if $\alpha \in
L^1(\mathbb {Z}^d)$ and $Z(\alpha)$ is $q$-thin, then $T^q[\alpha] =
L^q(G)$.  Here our $q$ is used in place of Edwards's $p$, and our $p$
is used in place of his $p'$.

We are now ready to prove Theorem \ref{Tabelian3}.  Suppose
$Z(\alpha)$ satisfies the hypothesis of the theorem. Let $\beta \in
L^p(\mathbb {Z}^d) \setminus 0$ such that $\alpha *\beta = 0$ and
$p \le \frac {2(d-\nu)}{d-1-\nu}$.  Since
$\frac {2(d-\nu)}{d-1-\nu} > 1$ and increasing $p$ retains the
property $\beta \in L^p(\mathbb {Z}^d)$, we may assume that $p>1$.
Then $\tilde{\alpha} * \tilde{\beta} = 0$
and using Proposition \ref{Pchar},
we see that $T^q[\alpha] \ne L^q(\mathbb {Z}^d)$.  But \cite[theorem
2.2]{edwards} tells us that $Z(\alpha)$ is not $q$-thin, and this
contradicts \cite[theorem 1]{guo}.

Conversely, let $T$ be a smooth, nonzero mass density on $Z(\alpha)$
vanishing near the boundary of $Z(\alpha)$.  Using \cite[theorem
3]{guo}, we can construct $\beta \in
L^p(\mathbb {R}^d) \setminus 0$ for $p > \frac{2(d-\nu)}{d-1-\nu}$
such that $\hat{\beta} = T$.  Then $\supp \hat{\beta} \subseteq
Z(\alpha)$, that is $\alpha \beta = 0$.
An application of Proposition \ref{Plinind} completes the proof of
Theorem \ref{Tabelian3}.

\section{Free groups and $p$-zero divisors}

Throughout this section, $2 \le k \in \mathbb {Z}$.  In
\cite{linnell} it was proven that if $0 \ne \alpha \in \mathbb
{C}F_k$, then $\alpha$ is not a 2-zero divisor.  In this section we
will give explicit examples to show that this result cannot be
extended to $L^p(F_k)$ for any $p>2$.  We will conclude this section
by giving sufficient conditions for elements of $L_r^1(F_k)$, the
radial functions of $L^1(F_k)$ as defined below,
to be $p$-zero divisors.

Any element $x$ of $F_k$ has a unique expression as a finite product
of generators and their inverses, which does not contain any two
adjacent factors $w w^{-1}$ or $w^{-1}w$.  The number of factors in
$x$ is called the \emph{length} of $x$ and is denoted by $|x|$.

A function in $L^{\infty} (F_k)$ will be called radial if its value
depends only on $|x|$.  Let $E_n = \{x \in F_k \mid |x| = n \}$, and
let $e_n$ indicate the cardinality of $E_n$.  Then $e_n =
2k(2k-1)^{n-1}$ for $n \ge 1$, and $e_0 = 1$.  Let $\chi_n$ denote
the characteristic function of $E_n$, so as an element of $\mathbb
{C}F_k$ we have $\chi_n = \sum_{|x|=n} x$.  Then every radial function
has the form $\sum_{n=0}^{\infty} a_n \chi_n$ where $a_n \in \mathbb
{C}$.  Let $L^p_r(F_k)$ denote the radial functions contained in
$L^p(F_k)$ and let $(\mathbb {C} F_k)_r$ denote the radial functions
contained in $\mathbb {C} F_k$.  Then $L^p_r(F_k)$ is the closure of
$(\mathbb {C} F_k)_r$ in $L^p(F_k)$.  Let $\omega = \sqrt{2k-1}$.
It was shown in \cite[chapter 3]{figapicardello} that
\begin{align*}
\chi_1 * \chi_1 &= \chi_2 + 2k * \chi_0 \\
\chi_1 * \chi_n &= \chi_{n+1} + \omega^2 \chi_{n-1}, \quad n\ge 2,
\end{align*}
hence $L^1_r(F_k)$ is a commutative algebra which is generated by
$\chi_0$ and $\chi_1$.

Later we will need the following elementary result.
\begin{Lem} \label{Ly}
Let $x,y \in F_k$ with $|x| = |y|$, and let $0 \le m,n \in \mathbb
{Z}$.  Then $\langle \chi_m * x, \chi_n \rangle = \langle \chi_m * y,
\chi_n \rangle$.
\end{Lem}
\begin{proof}
We have $\langle \chi_m * x, \chi_n \rangle
= \langle x, \chi_m^* * \chi_n \rangle
= \langle x, \chi_m * \chi_n \rangle$.
By the above remarks, $\chi_m * \chi_n$ is a sum
of elements of the form $\chi_r$.  Therefore we need only prove that
$\langle x, \chi_r \rangle = \langle y, \chi_r \rangle$.  But
$\langle x, \chi_r \rangle = 1$ if $m=r$ and 0 if $m \ne r$, and the
result follows.
\end{proof}

Let $\alpha$ be a complex-valued function on $F_k$.  Set
\[
a_n (\alpha) = \frac{1}{e_n} \sum_{x \in E_n} \alpha (x)
\]
and denote by $P(\alpha)$ the radial function $\sum_{n=0}^{\infty}
a_n(\alpha) \chi_n$.

\begin{Lem} \label{Lradial1}
Let $1 \le p \in \mathbb {R}$ or $p = \infty$, let
$\alpha \in L^1_r (F_k)$, and let $\beta \in L^p(F_k)$.  If
$\alpha * \beta = 0$, then $\alpha * P(\beta) = 0$.
\end{Lem}
\begin{proof}
Let $f,h \in \mathbb {C}F_k$.  It was shown in \cite[lemma
6.1]{pytlik} that $ P(f) * P(h) = P(P(f)*h) $.
Write $\beta = \sum_{g \in F_k} b_g g$.  If $p \ne \infty$ and
$0 \le a_1, \dots, a_n \in
\mathbb {R}$, then by Jensen's inequality \cite[p.~189]{jensen}
applied to the function $x^p$ for $x > 0$,
\[
\left(\frac{a_1 + \dots + a_n}{n}\right)^p \le
\frac{a_1^p + \dots + a_n^p}{n},
\]
consequently
\[
\|P(\beta)\|_p^p =
\sum_{n=0}^\infty e_n \left| \frac{1}{e_n} \sum_{|g| = n} b_g
\right|^p \le \sum_{g \in F_k} |b_g|^p = \| \beta \|_p^p.
\]
Therefore $P$ is a continuous map from $L^p(F_k)$ into $L^p_r(F_k)$
for $p \ne \infty$.  It is also continuous for $p = \infty$.
The lemma follows because the map $L^1(G) \times L^p(G) \to L^p(G)$
is continuous; specifically $\| \alpha * \beta \|_p
\le \| \alpha \|_1 \|\beta \|_p$.
\end{proof}

For $n \in \mathbb {Z}_{\ge 0}$, define polynomials $P_n$ by
\begin{align*}
&P_0(z) = 1,\quad P_1(z) = z,\quad P_2(z) = z^2 -2k \\
\text{and }
&P_{n+1}(z) = zP_n(z) - \omega^2 P_{n-1}(z) \text{ for } n\ge 2.
\end{align*}
Let $\alpha = \sum_{n=0}^{\infty} a_n \chi_n \in L^1_r (F_k)$.  In
\cite{pytlik}, Pytlik shows the following.
\begin{enumerate}
\item
$X = \{x + iy \in \mathbb {C} \mid (\frac {x}{2k})^2 + (\frac
{y}{2k-2})^2 \le 1 \}$ is the spectrum of $L^1_r(F_k)$.

\item The Gelfand transform of $\alpha$ is given by
$\hat{\alpha}(z) =
\sum_{n=0}^{\infty} a_n P_n(z)$ for $z \in X$.
\end{enumerate}
Let $Z(\alpha) = \{z \in X \mid \hat{\alpha} (z) = 0 \}$.
For $z \in X$ we define $\phi_z \in L_r^{\infty}(F_k)$, the space of
continuous linear functionals on $L_r^1(F_k)$
\cite[p.~34]{cohen}, by
\[
\phi_z = \sum_{n=0}^{\infty} \frac{P_n(z)}{e_n} \chi_n.
\]

We can now state
\begin{Lem} \label{Lradial2}
Let $\alpha \in L^1_r(F_k)$ and let $z \in X$.  Then $\alpha
* \overline{\phi_z} = 0$ if and only if $z \in Z(\alpha)$.
\end{Lem}
\begin{proof}
Let $\beta \in L^1_r(F_k)$ and write $\beta = \sum_{m=0}^{\infty} b_m
\chi_m$.  Then 
\begin{align*}
\langle \beta, \overline{\phi_z} \rangle &= \sum_{m,n}
\frac{ b_m P_n(z)}{e_n}
\langle \chi_m, \chi_n \rangle \\
&= \sum_n b_n  P_n(z) =  \hat{\beta} (z).
\end{align*}
Applying this in the case $\beta = \alpha * \chi_n$, we obtain
$\langle \alpha * \chi_n,
\overline{\phi_z} \rangle = \hat{\alpha}
(z)  P_n (z)$.  Using Lemma \ref{Ly}, we deduce that
if $y \in F_k$
and $|y| = n$, then $\langle \alpha * y , \phi_z \rangle
=  \hat{\alpha} (z)
P_n(z)/e_n$.  Since $\alpha = \tilde{\alpha}$, the
result now follows from Lemma \ref{Lperp}.
\end{proof}

If $\alpha \in L^1_r(F_k)$, we shall say that $\alpha *\chi_n$ is a
radial translate of $\alpha$.  We then set $TR^1[\alpha]$ equal
to the closure in $L^1_r(F_k)$ of the set of linear combinations of
radial translates of $\alpha$.

\begin{Prop} \label{Pradial}
Let $\alpha \in L^1_r (F_k)$.  Then $\alpha *\beta \ne 0$ for
all $\beta \in L^{\infty}(F_k) \setminus 0$
if and only if $Z(\alpha) = \emptyset$.
\end{Prop}
\begin{proof}
If $z \in Z(\alpha)$, then $\phi_z \in L^{\infty}(F_k)
\setminus 0$ and
$\alpha * \overline{\phi_z} = 0$ by Lemma \ref{Lradial2}.

Conversely suppose there exists $\beta
\in L^{\infty}(F_k) \setminus 0$ such
that $\alpha * \beta = 0$.  Then $\beta_y \ne 0$ for some $y \in
F_k$, so replacing $\beta$ with $\beta
* y^{-1}$, we may assume that
$P(\beta) \ne 0$.  If $\gamma = \bar {\beta}$, then $\alpha *
\bar {\gamma} = 0$ and $P(\gamma) \ne 0$.
Using Lemma \ref{Lradial1} we see that $\alpha *
\overline{P(\gamma)} = 0$, and we deduce from Lemma
\ref{Lperp} that $\langle \alpha_y, P(\gamma) \rangle = 0$ for all $y
\in F_k$.  It follows that $\langle \alpha * \chi_n, P(\gamma)
\rangle = 0$ for all $n \in \mathbb {Z}_{\ge 0}$, consequently
$TR^1[\alpha] \ne L^1_r(F_k)$.  Let $J$ be
a maximal ideal in $L^1_r(F_k)$ which contains $TR^1[\alpha]$.  By
Gelfand theory there exists $z \in X$ such that $J = \{\delta \in
L^1_r(F_k) \mid \hat{\delta}(z) = 0 \}$, so $z \in Z(\gamma)$.
\end{proof}

We can now state
\begin{Ex} \label{Ezero1}
Let $k \ge 2$.  Then $\chi_1$ is a $p$-zero divisor for all $p>2$.
\end{Ex}
\begin{proof}
Since $0 \in Z(\chi_1)$, we see from Lemma \ref{Lradial2}
that $\chi_1 * \phi_0 = 0$.  Of
course $\phi_0 \ne 0$.  We now prove the stronger statement that
$\phi_0 \in L^p(F_k)$ for all $p>2$.  We have
\[
\phi_0 = \sum_{n=0}^{\infty} \frac{P_n(0)}{e_n} \chi_n =
\sum_{n=0}^{\infty} \frac{(-1)^n}{(2k-1)^n} \chi_{2n}.
\]
Therefore
\begin{align*}
\sum_{g \in F_k} |\phi_0(g)|^p &=
1 + \sum_{n=1}^{\infty} \frac{e_{2n}}{(2k-1)^{pn}}
= 1 + \sum_{n=1}^{\infty} \frac{2k(2k-1)^{2n-1}}{(2k-1)^{pn}} \\
&= 1 + \frac{2k}{2k-1} \sum_{n=1}^{\infty} \frac{1}{(2k-1)^{n(p-2)}}
\end{align*}
and the result follows.
\end{proof}

We can use the above result to
prove that the nonsymmetric sum of generators in $F_k$ is
a $p$-zero divisor for all $p>2$ in the case $k$ is even
and $k>2$.  Specifically we have
\begin{Ex} \label{Ezero2}
Let $k>3$ and let $\{ x_1, \dots , x_k \}$ be a set of generators
for $F_k$.  If $k$ is even,
then $x_1 + \dots + x_k$ is a $p$-zero divisor for
all $p>2$.
\end{Ex}

To establish this, we need some results about free groups.
\begin{Lem} \label{Lfree}
Let $0 < n \in \mathbb {Z}$ and let $F$ be the free group on $x_1,
\dots, x_n $.  Then no nontrivial word in the $2n-1$ elements
$x_1^2, \dots, x_n^2, x_1x_2, x_2x_3, \dots, x_{n-1}x_n$ is the
identity; in particular these $2n-1$ elements generate a free group
of rank $2n-1$.
\end{Lem}
\begin{proof}
We shall use induction on $n$, so assume that the result is true with
$n-1$ in place of $n$.
Let $T$ denote the Cayley graph of $F$ with respect to the generators
$x_1, \dots, x_n$.  Thus the vertices of $T$ are the elements of $F$,
and $f,g \in F$ are joined by an edge if and only if $f = gx_i^{\pm
1}$ for some $i$.  Also $F$ acts by left multiplication on $T$.
Suppose a nontrivial word in $x_1^2, \dots, x_n^2, x_1x_2, x_2x_3,
\dots, x_{n-1}x_n$ is the identity, and choose such a word $w$ with
shortest possible length.  

Note that $w$ must involve $x_n^2$, because $F$ is the free product
of the group generated by $x_1, \dots, x_{n-1}$ and the group
generated by $x_{n-1}x_n$.  By conjugating and taking inverses if
necessary, we may assume
without loss of generality that $w$ ends with $x_n^2$.  

Write $w = w_m \dots w_1$, where $w_1 = x_n^2$, and each of
the $w_i$ are one of the above $2n-1$ elements.  Let us consider the
path whose $(2i+1)$th vertex is $w_i \dots w_1 1$.  Note that
$w1 = 1$, but $w_i \dots w_11 \ne 1$ for $0 < i < m$.

Observe that the path of length 2 from $x_n^2$ to
$w_2 x_n^2$ cannot go through $x_n$ (just go through the $2n-1$
possibilities for $w_2$, noting that $w_2 \ne x_n^{-2}$).  Now remove
the edge joining $x_n$ and $x_n^2$.  Since $T$ is a tree \cite[I.8.2
theorem]{dicks}, the resulting graph will become two trees; one
component $T_1$ containing 1 and the other component $T_2$ containing
$x_n^2$.  Since the length 2 path from $x_n^2$ to $w_2 x_n^2$ did not
go though $x_n$, for $i \ge 1$ the path $w_i \dots w_2(w_1 1)$
remains in $T_2$ at least
until it passes through $x_n^2$ again.
Also the path must pass through $x_n^2$ again in order to get back to
1.  Since the paths $w_i \dots (w_11)$ all have even length
(all the $w_i$ are words of length 2), it follows that
$w_l \dots w_1 1 = x_n^2$ for some
$l \in \mathbb {Z}$, where $2 \le l < m$.
We deduce that $w_l \dots w_2 = 1$, which contradicts the
minimality of the length of $w$.
\end{proof}

\begin{Cor} \label{Cfree1}
Let $n \in \mathbb {Z}_{\ge 1}$ and let $F$ be the free group on
$x_1, \dots , x_n$.  Then no nontrivial word in the $2n-1$ elements
$x_1^2 , \dots, x_n^2, x_1^{-1}x_2, x_2^{-1}x_3, \dots, x_{n-1}^{-1}
x_n$ is the identity; in particular these $2n-1$ elements generate a
free group of rank $2n -1$.
\end{Cor}
\begin{proof}
This follows immediately from Lemma \ref{Lfree}: replace $x_i
x_{i+1}$ with $x_i^{-2} x_i x_{i+1}$ for all $i < n$.
\end{proof}

\begin{Cor} \label{Cfree2}
Let $n \in \mathbb {Z}_{\ge 1}$ and let $F$ be the free group on
$x_1, \dots, x_n, w$.  Then the elements $wx_1, wx_1^{-1}, \dots,
wx_n, wx_n^{-1}$ generate a free subgroup of rank $2n$.
\end{Cor}
\begin{proof}
The above elements generate the subgroup generated by
\[
x_1^2, \dots, x_n^2, x_1^{-1}x_2, x_2^{-1}x_3, \dots,
x_{n-1}^{-1}x_n, wx_1 .
\]
The result follows from Corollary
\ref{Cfree1}.
\end{proof}

\begin{proof}[Proof of Example \ref{Ezero2}]
Let $G = F_k$ and let
$F$ be the free group on $y_1, \dots , y_k, w$.  By Corollary
\ref{Cfree2} there is a monomorphism $\theta \colon G \to F$
determined by the formula
\[
\theta (x_1) = wy_1, \quad
\theta (x_2) = wy_1^{-1}, \quad
\dots, \quad
\theta (x_k) = wy_{k/2}^{-1}.
\]
Note that $\theta$ induces a Banach space monomorphism $L^p(G) \to
L^p(F)$.
Set $\alpha = wy_1 + wy_1^{-1} + \dots + wy_{k/2} + wy_{k/2}^{-1}$.
Since $y_1 + y_1^{-1} + \dots + y_{k/2} + y_{k/2}^{-1}$ is a $p$-zero
divisor by Example \ref{Ezero1}, we see that $\alpha$
is a $p$-zero divisor, say $\alpha * \beta = 0$ where $0 \ne \beta \in
L^p(F)$.  Write $F = \bigcup_{t \in T}
\theta (G)t$ where $T$ is a right
transversal for $\theta (G)$ in $F$.
Then $\beta = \sum_{t\in T} \beta_t t$
where $\beta_t \in L^p(\theta (G))$ for all $t$.  Also $\alpha *
\beta_t = 0$ for all $t$ and $\beta_s \ne 0$ for some $s \in T$.
Define $\gamma \in L^p(G)$ by $\theta (\gamma) = \beta_s$.  Then
$0 \ne \gamma \in L^p(G)$ and $(x_1 + \dots + x_k)
* \gamma = 0$ as required.
\end{proof}

We conclude with some information on the existence of $p$-zero
divisors in $L^1_r(F_k)$.  Let $\alpha \in L^1_r(F_k)$ and define
$p(\alpha)$ as follows.  If $Z(\alpha) \cap (-2k,2k) = \emptyset$,
then set $p(\alpha) = \infty$.  If $Z(\alpha) \cap (-2k,2k) \ne
\emptyset$, then set $m(\alpha) = \min \{|t| \mid t \in Z(\alpha)
\cap (-2k, 2k) \}$.  If $m(\alpha) \in [0, 2\omega]$, then set
$p(\alpha) = 2$.  Finally if $m(\alpha) \in (2\omega, 2k)$, then let
$p(\alpha)$ be the positive root of the following equation in $p$:
\[
m(\alpha) = \sqrt{2k-1} \bigl( (2k-1)^{\frac{1}{2} - \frac{1}{p}} +
(2k-1) ^{\frac{1}{p} - \frac{1}{2}} \bigr).
\]
We can now state
\begin{Prop}
Let $\alpha \in L^1_r (F_k)$.  Then $\alpha $ is a $p$-zero divisor
for all $p > p(\alpha)$.
\end{Prop}
\begin{proof}
Let $t \in (-2k, 2k)$
such that $m(\alpha) = |t|$ and suppose $p >  p(\alpha)$.
Since $\phi_t$ is a positive definite
function by \cite[lemma 6.1]{pytlik}, we can apply \cite[theorem
2(a)]{cohen} to deduce that $\phi_t \in L^p_r(F_k)$.  By Lemma
\ref{Lradial2} $\alpha * \phi_t = 0$ and the result is proven.
\end{proof}

\bibliographystyle{plain}

\end{document}